\documentclass[11pt]{amsart}
\usepackage{}
\usepackage{amssymb}
\theoremstyle{plain}
\newtheorem{Thm}{Theorem}

\newtheorem{Def}[Thm]{Definition}

\begin{document}

\title[injectivity radius bound]
{injectivity radius bound and good leaf foliation in 3-manifolds}

\author{Li Ma}

\address{ Department of mathematics \\
Henan Normal university \\
Xinxiang, 453007 \\
China}

\email{lma@tsinghua.edu.cn}

\thanks{ The research is partially supported by the National Natural Science
Foundation of China (N0.11271111)}

\begin{abstract}
In this short note, we study the injectivity radius bound for three
dimensional complete and non-compact Riemannian manifold with good
leaf foliations and with bounded curvature up to first order. We
obtain the injectivity bound by using the minimal surface theory and
the Gauss-Bonnet theorem.


\end{abstract}

 \maketitle

\section{ main result and its proof}

In order to use Hamilton's Cheeger-Gromov convergence theorem of
Ricci flows, we need to know the uniform injectivity radius lower
bound of the underlying Riemannian metrics. That is to say, we need
to handle the local collapsing of a family of Riemannian manifolds.
Generally speaking, this is a hard work. However, there are two most
useful approaches to this goal. One is to find the volume lower
bound via the use of results in the paper of Cheeger-Gromov-Taylor
in 1982 \cite{C}. The other is to use Perelman's W-functional
\cite{P02} on Riemannian manifolds.

  In this note, we propose the minimal surface method on three dimensional complete non-compact
  Riemannian manifold $(M,g)$ with bounded curvature up to first order and with a "good leaf foliation".
  Here, the bounded curvature up to first
  order just means that both curvature and its covariant derivatives of the curvature are
  bounded. The concept "good leaf foliation" is defined below.

  \begin{Def}\label{foliation}
  We say $M$ has a "good leaf foliation" if there is a proper function $f:M\to [0,\infty)$  such that $\Omega_a=\{f\leq a\}$
  expanding $M$ as $a\to\infty$ and for almost all $a\geq 0$, the leaf
  $\{f(x)=a\}$ has positive mean curvature (or the so called mean
  curvature  convexity ).
\end{Def}

  Such a foliation may be called a "mean curvature convexity foliation" since it also occurs in the study of mean curvature
  flow. In history, mean curvature convexity (or the convexity) plays an important role in geometry \cite{B}.
  See also the recent work of M.Gromov \cite{G}.

The below is our main result.
\begin{Thm}\label{mali} Let $(M,g)$ be a three dimensional simply connected complete non-compact
  Riemannian manifold with bounded curvature up to first order and with a "good leaf
  foliation". Then there a uniform constant $b>0$ such that $injrad
  (M,g)\geq b$.
\end{Thm}
 The proof is to use the minimal surface theory on three dimensional Riemannian manifolds
 and the main argument is along that of in \cite{MZ} but simpler.

\begin{proof} (of Theorem \ref{mali}).
 We argue by
contradiction.

 Assume that we have uniform lower bound of the injectivity radius on $(M,g)$.
 That is to say, there are $r_j:=injrad(x_j,g)\to 0$ for some sequence $(x_j)\subset M$ going
 to infinity. This implies that we can find for each $x_j$ a
 geodesic loop $\gamma_j$ which may have corner at $x_j$ and $L(\gamma_j, g)\to 0$. We can
 smooth $\gamma_j$ at $x_j$ and get a approximated geodesic loop,
 which is still denoted by $\gamma_j$. Then we can span a
 minimizing disk $D_j$ over the curve $\gamma_j$ (which is the so called Douglas-Morrey minimal surface). The minimal
 disk $D_j$ can not touch any point in the level set $f(x)=0$
 since it has the property that $\Delta f>0$ on $f(x)=0$.

 We claim that the minimal disk is between two level sets of $f$ of uniform width,
 saying the bounded region $\Omega_j$, which is exact the region obtained by two leafs bounding the circle $\gamma_j$.
  We can say that the region $\Omega_j$ is the region between $f=a_j+L_j$ and $f=a_j$ with $L_j\leq
 \frac{1}{2}L(\gamma_j,g)\leq 3r_j$. The reason is below. If we have an interior point
 $p_j\in D_j$ such that $f(p_j)>a_j+L_j$, then setting $m_j:= \sup\{c\in f(D_j)\}$,
 we can use some level set $f(x)=m_j-\epsilon:=t_j>a_j+L_j$ to cut from the disk $D_j$ a small
 loop $\tau_j$ on the leaf $f(x)=t_j$.
 Using Theorem 1 in the paper of Meeks-Yau \cite{MY}, we can construct another minimizing disk $B_j$
 spanning $\tau_j$ such that $B_j\subset \{f(x)\leq t_j\}$.
 Note that for $\epsilon>0$ small, $D_j \bigcap \{f(x)\geq t_j\}$ is a minimal disk. That is to say, we obtain two
 different minimizing disks spanning the small loop $\tau_j$, which contradicts
 to the uniqueness result of minimal surface for small enough loops.

 We now consider the rescaling metric $g_j:=\frac{1}{r_j^2}g$  of the metric $g$ at
 $x_j$. Recall that $|Rm(g_j)|_{g_j}=r_j^2|Rm(g)|_g\to 0$ as
 $j\to\infty$.
 Then we know by our curvature assumption that the Gromov-Hausdorff limit $(M_\infty, g_\infty, x_\infty)$ of
 the pointed Riemannian manifolds sequence $(M,g_j,x_j)$ is one of the flat space $R^\times T^2$, and $R^2\times
 S^1$. Note that $M_\infty=R^3$ is impossible since $injrad(x_\infty, g_\infty)=1$. Now
 the rescaling minimal disk $\hat{D}_j$ of the minimal disk $D_j$ is bounded
 in a region of maximum diameter no more than two in the metric
 $g_j$. Hence its limit $D_\infty$ in $(M_\infty,g_\infty,x_\infty)$
 is a flat disk. By this we know that the area $A(\hat{D}_j)$ is uniformly
 bounded. Then we can use the argument in \cite{MZ} to get a
 contradiction. For completeness, we repeat it here.

 In below we rewrite $D_j$ as $\hat{D}_j$ for simplicity. On the surface
$D_j$ with the unit normal vector $N$ and the second fundamental
form $A_j$ of $D_j\subset M$, we know that
$$
Rc_{g_j}(N,N)+\frac{1}{2}|A_j|^2=\frac{R(g_j)}{2}-K_j.
$$
Here $K_j$ is the Gauss curvature of the minimal surface $D_j$ and
$Rc_{g_j}$ and $R(g_j)$ are the Ricci and scalar curvatures of the
metric $g_j$ respectively. Then using $Rc_{g_j}\to 0$ we have
$$
\frac{1}{2}|A|^2+K_j\leq\frac{R(g_j)}{2}-Rc_{g_j}(N,N)\to 0
$$
and by using the uniform area bound of $D_j$,
$$
\int_{D_j}K_j\leq \frac{1}{2}\int_{D_j}R(g_j)\leq \sup_{D_j}
|R(g_j)-2Rc_{g_j}(N,N)| A(D_j)\to 0.
$$

Note that $\int_{\tau_j}k_{g_j}\to 0$ where $k(g_j)$ is the geodesic
curvature of the almost geodesic loop $\tau_j$.

 Recall that the Gauss-Bonnet formula (and here we know that $D_j$'s are embedded)
 $$
\int_{D_j}K_j+\int_{\tau_j}k_{g_j}=2\pi,
 $$
 which gives us a contradiction since both terms in left side tend
 to zero or less than  zero as $j\to\infty$.

 This completes the proof of Theorem \ref{mali}.

\end{proof}

\textbf{Remark}. One important example which has the leaf foliation
as in  the conditions of Theorem \ref{mali} is the three dimensional
complete non-compact riemannian manifold with bounded and positive
Ricci curvature as used in \cite{MZ}. In the latter case, one simply
takes $f$ as the Busemann function defined by a ray.

{\bf Acknowledgement}. This work is done while the author is
visiting IHES, France and the authors would like to thank the
hospitality of IHES.


\begin{thebibliography}{10}

\bibitem{B} M.Berger, \emph{Geometry revealed}, Springer, 2010

\bibitem{C}
 B.Chow, P.Lu, L.Ni. \emph{Hamilton's Ricci Flow}. Beijing: Science
Press, American Mathematical Society, 2006

\bibitem{G} M.Gromov, \emph{Plateau-Stein manifolds}, IHES preprint, 2013.

\bibitem{MZ} Li Ma, A.Q.Zhu, \emph{Injectivity radius bound of Ricci flow with
positive Ricci curvature and applications}, Front. Math. China 2013,
8(5): 1129-1137, DOI 10.1007/s11464-013-0296-8

\bibitem{MY} W.W.Meeks II, S.T.Yau, \emph{The existence of embedded minimal
surfaces and the problem of uniqueness}, Math.Z., 179(1982)151-168

\bibitem{P02} Grisha Perelman,
\emph{The entropy formula for the Ricci flow and its geometric
applications}, http://arxiv.org/abs/math/0211159v1

\end{thebibliography}
\end{document}